\documentclass[12pt,twoside]{article}
\usepackage{amsmath,amsbsy,amsfonts}
\newtheorem{theorem}{Theorem}[section]
\newtheorem{lemma}[theorem]{Lemma}
\newtheorem{proposition}[theorem]{Proposition}

\newtheorem{remark}[theorem]{Remark}

\newcommand{\qed}{\hfill\rule{2mm}{2mm}}

\title{Existence of solutions for a class of $p(x)$-laplacian equations involving a concave-convex nonlinearity with critical growth in $\mathbb{R}^{N}$\footnote{Partially supported by INCT-MAT and PROCAD}}

\author{Claudianor O. Alves\footnote{C.O. Alves was partially supported by CNPq/Brazil
303080/2009-4, e-mail:coalves@dme.ufcg.edu.br}
\,\,\, and \,\,\,  Marcelo C. Ferreira \footnote{e-mail:marcelo@dme.ufcg.edu.br}\\
Universidade Federal de Campina Grande\\
Unidade Acad\^emica de Matem\'atica\\
CEP:58429-900, Campina Grande - PB, Brazil.}

\date{}

\begin{document}
\maketitle

{\scriptsize{\bf 2000 Mathematics Subject Classification:} 35A15, 35H30, 35B33.}

{\scriptsize{\bf Keywords:} Variational Methods, $p(x)$-laplacian, Critical Growth }

\begin{abstract}

We prove the existence of solutions for a class of quasilinear problems involving variable exponents and with nonlinearity having critical growth. The main tool used is the variational method, more precisely, Ekeland's Variational Principle and the Mountain Pass Theorem.

\end{abstract}

\section{Introduction}

The present paper concerns with the existence of solutions for the following class of quasilinear problems involving variable exponents
   $$
        \begin{cases}
            -\Delta_{p(x)} u + V(x) u^{ p(x)-1 } = \lambda h(x) u^{ r(x)-1 } + \mu u^{ q(x)-1 } + u^{ p^{ \ast }(x)-1 }, \, \mathbb{R}^{N} \\
            u \ge 0 \,\,\, \mbox{and} \,\, u \not=0, \,\, \mathbb{R}^{N}\\
            u \in W^{ 1,p(x) } \big( \mathbb{R}^{N} \big),
         \end{cases} \eqno{(P)}
   $$
where $\Delta_{p(x)}$ is the $p(x)$-laplacian operator given by
   $$
      \Delta_{p(x)} u= \text{div}(|\nabla u|^{p(x)-2} \nabla u),
   $$
   $ \lambda , \mu  $ are positive parameters, $ p \colon \mathbb{R}^{N} \to \mathbb{R} $ is a Lipschitz continuous function, $ V, q, r :  \mathbb{R}^{N} \to \mathbb{R} $ are continuous functions  and  $ h $ is a nonnegative function in $ L^{ \Theta(x) } \big( \mathbb{R}^{N} \big) $ with
   $$
       \Theta(x) = \dfrac{Np(x)}{Np(x)-r(x) \big( N-p(x) \big)}.
   $$
   \noindent Moreover, the functions $p,q$ and $V$ are $\mathbb{Z}^{N}$-periodic, that is
   $$
     p(x+y)=p(x), \, q(x+y)=q(x), \, V(x+y)=V(x) \,\,\, \forall x \in \mathbb{R}^{N} \,\,\, \mbox{and} \,\,\, \forall y \in \mathbb{Z}^{N}
     \eqno{(H_0)}
  $$
and we also assume that
  $$
     1 < p_- \leq p(x) \leq p_+ < N \,\,\, \forall x \in \mathbb{R}^{N}. \eqno{(H_1)}
  $$
  $$
      1 < r_- \le r_+ < p_- \le p_+ < q_- \leq  q(x) \ll p^{ \ast }(x), \, \forall x \in \mathbb{R}^{N}. \eqno{(H_2)}
  $$
  $$
     \inf_{ x \in \mathbb{R}^{N} } V(x) = V_0>0.  \eqno{(V_0)}
  $$
  Here, the notation $u(x)  \ll  v(x)$  means that  $\displaystyle \inf_{x\in\mathbb{R}^{N}}(u(x)-v(x))>0$, \linebreak $u_- = \displaystyle \mbox{ess}\inf_{x \in \mathbb{R}^{N}}u(x)$, $u_+ = \displaystyle \mbox{ess}\sup_{x \in \mathbb{R}^{N}}u(x)$ and $u^{ \ast }(x)=\frac{Nu(x)}{N-u(x)} \, \, \, \forall x \in \mathbb R^N$.

  \vspace{0.5 cm}

  Partial Differential Equations involving the $ p(x) $-laplacian arise, for instance, as a mathematical model for problems involving electrorheological fluids and image restorations,
  see \cite{Acerbi1,Acerbi2,Antontsev,CLions,Chen,Ru}. This explains the intense research on this subject in the last decades. Regarding to the application of variational methods
  in order to solve $ p(x) $-laplacian problems, many research were already done when the nonlinearities have a subcritical growth, see for example,
 \cite{AlvesSouto,Alves2,CP,FanZa,FH,FZ1,Fan1,Fan,MR} and references therein. However, when the growth involves some criticality, some articles just began appear recently, see the papers
 due to Alves \& Souto \cite{AlvesSouto},  Alves \cite{Alves3}, Alves \& Ferreira \cite{AlvesFerreira}, Bonder \& Silva  \cite{BS}, Bonder, Saintier and Silva \cite{BSS,BSS1}, Fu \& Zhang \cite{FuZa,FuZa1}, Shang \& Wang \cite{SW} and references therein.

  In  \cite{Alves}, Alves has studied the existence of solutions for the following class of quasilinear problems:
  $$
        \begin{cases}
        -\Delta_p u = \lambda g(x) u^{ r-1 } + u^{ p^*-1 }, \, \mathbb R^N \\
        u \geq 0, \, u \neq 0 \\
        u \in D^{ 1,p }( \mathbb R^N ),
     \end{cases} \eqno{(P_0) }
   $$
   where $\lambda >0$,  $ 2 \leq p \leq N $,  $ 1 < r < p $ and $ g $ is a nonnegative function belonging to $ L^{ \theta }( \mathbb R^N ) $ with
   $$
       \theta = \frac{Np}{Np-r(N-p)}.
   $$
In \cite{Alves}, by using variational methods, more precisely, Mountain Pass Theorem and Ekeland's Variational Principle, the existence of two solutions has been established when $\lambda$ is small enough.  In the literature, we can find a lot of papers related to problem $(P_0)$ involving bounded or unbounded domains, see for example, \cite{ ABCerami,BarWil, CBZ,GonAlves,GV,Pan,Tarantello,W}. However, involving variable exponents, the authors know only the paper \cite{BS}, where the nonlinearity has a behavior like concave-convex and the domain is bounded.

Motivated by the above informations, we prove that similar results to that found in \cite{Alves} also hold for the case where the exponents are variable. More precisely, we have showing that the energy functional $I:W^{ 1,p(x) }( \mathbb{R}^{N} ) \to \mathbb{R} $ associated with $ (P) $, which is given by
  \begin{multline*}
     I( u ) = \int_{\mathbb{R}^{N}} \frac{1}{p(x)} \left( | \nabla u |^{ p(x) } + V(x) | u |^{ p(x) } \right) - \lambda \int_{\mathbb{R}^{N}} \frac{h(x)}{r(x) } ( u^+ )^{ r(x) } \\
      - \mu \int_{\mathbb{R}^{N}} \frac{1}{q(x)} ( u^+ )^{ q(x) } - \int_{\mathbb{R}^{N}} \frac{1}{p^{ \ast }(x) } ( u^+ )^{ p^{ \ast}(x) },
  \end{multline*}
has two critical points for each $\mu $ large enough and $\lambda$ small enough.

\vspace{0.5 cm}

   Our main theorem is the following

\begin{theorem}  \label{T1}
   There exists $ \mu^{ \star } > 0 $ such that for each $ \mu \geq \mu^{ \star } $, there is $ \lambda_\mu = \lambda (\mu) > 0 $ such that problem $ (P) $ has two solutions
   $ \Psi_1, \Psi_2 \in W^{ 1,p(x) } ( \mathbb R^N ) $ with
   $$
       I(\Psi_2)<0<I(\Psi_1),
   $$
   for all $ \lambda \in (0, \lambda_\mu)$. \\
\end{theorem}

The Theorem \ref{T1} is an immediate consequence of Theorems \ref{T2} and \ref{T3}, which were proved in Sections 4 and 5 respectively. In the proof of the above results, we have used a result found in \cite{AlvesFerreira},  which shows that if $(H_0)-(H_2)$ hold, the problem
     $$
          \begin{cases}
       -\Delta_{p(x)}u + V(x)| u |^{ p(x)-2 } u = \mu |u|^{ q(x)-2 } u + | u |^{ p^{ \ast }(x)-2 } u,  \, \mathbb{R}^N \\
       u \neq 0 \,\,\, \mbox{and} \,\,\, u \in W^{ 1,p(x)} \big( \mathbb{R}^{N} \big),
     \end{cases} \eqno{(P_\mu) }
  $$
 has a {\it ground state solution}, that is, the mountain pass level of the energy functional associated with $(P_\mu)$ is a critical value.

 We recall that the energy functional $I_\mu: W^{1,p(x)}(\mathbb{R}^{N}) \to \mathbb{R}$ associated to $(P_\mu)$ is given by
  $$
     I_\mu(u)=\int_{\mathbb{R}^{N}} \frac{1}{p(x)} \left( | \nabla u |^{ p(x) } + V(x) | u |^{ p(x) } \right) - \mu \int_{\mathbb{R}^{N}} \frac{1}{q(x)} | u |^{ q(x)} - \int_{\mathbb{R}^{N}} \frac{1}{p^{ \ast }(x)} |u|^{ p^{ \ast}(x) }.
  $$
Thus, if $c_\mu$ denotes the mountain pass level of $I_\mu$, we say that $\Psi \in W^{1,p(x)}(\mathbb{R}^{N})$ is a {\it ground state solution} of $(P_\mu)$ if
$$
I_\mu '(\Psi)=0 \,\,\, \mbox{and} \,\,\, I_\mu(\Psi)=c_\mu.
$$
In \cite{Alves3}, the below limit has been proved

   \begin{equation} \label{LIMITEDONIVEL}
      c_\mu \to 0, \,\,\, \mbox{as} \,\,\, \mu \to +\infty.
  \end{equation}
The above limit is a key point in our arguments, because in the present paper,  we will denote by $ \mu_0 > 0 $ a number such that
   \begin{equation} \label{ESTIMATIVAPASSO}
      c_\mu < \min \left \{ \gamma \left( \frac{1}{K} \right)^{ \frac{1}{\gamma} }, \frac{1}{2K^{ p_+ }} \nu \right \} \,\,\, \forall \mu \ge \mu_0,
  \end{equation}
  where
  \begin{equation} \label{theta}
     \gamma = 1/p_+ - 1/p^{ \ast}_-, \ \nu = 1/p_+ - 1/q_-,
  \end{equation}
  and $K \ge 1$ is fixed satisfying
  $$
    | u |_{ p^{ \ast }(x) } \le K \| u \|, \,\,\,\,  \forall u \in W^{ 1,p(x) } \big( \mathbb{R}^{N} \big).
  $$
Furthermore, standard arguments work to prove that the ground state solution $\Psi$ of $ (P_\mu) $ can be chosen nonnegative.

\vspace{0.5 cm}

\noindent {\bf Notation:} The following notations will be used in the present work: \\

\noindent $\bullet$ \,\, $C$ and $C_i$ will denote generic positive constant, which may vary from line to line.  \\

\noindent $\bullet$ \,\, In all the integrals we omit the symbol $dx$. \\

\noindent $\bullet$ \,\, $u^+(x) = \max\{ u(x), 0 \}$ \,\,\, and \,\,\, $u^-(x) = \min\{ u(x), 0 \}. $

\section{Variable exponent Lebesgue and Sobolev spaces}

    In this section, we recall some results on variable exponent Lebesgue and Sobolev spaces found in \cite{FSZ,FZ} and their references.

    Let $ z\in L^{\infty}(\mathbb{R}^{N}) $ with $ z_- \ge 1$. The variable exponent Lebesgue space $L^{z(x)}(\mathbb{R}^{N})$ is defined by
    \[
        L^{z(x)}(\mathbb{R}^{N})=\left\{ u:\mathbb{R}^{N} \to  \mathbb{R} \left\vert \,u\text{ is measurable and }\int_{\mathbb{R}^{N}}\left\vert u\right\vert^{z(x)}<\infty\right.\right\},
    \]
    endowed with the norm
    \[
       \left\vert u\right\vert _{z(x)}=\inf\left\{  \lambda>0\left\vert \,\int_{\mathbb{R}^{N}}\left\vert \frac{u}{\lambda}\right\vert ^{z(x)} \leq 1\right.  \right\}  \text{.}
    \]
    The variable exponent Sobolev space is defined by
    \[
        W^{1,z(x)}(\mathbb{R}^{N})=\left\{  u\in L^{z(x)}(\mathbb{R}^{N})\left\vert \,\left\vert \nabla u\right\vert \in L^{z(x)}(\mathbb{R}^{N})\right. \right\},
    \]
    with the norm
   \[
      \left\Vert u\right\Vert _{1,z(x)}=\left\vert u\right\vert_{z(x)}+\left\vert \nabla u\right\vert _{z(x)}\text{.}
   \]
   If $M \in L^{\infty}(\mathbb{R}^{N})$ satisfies $ M_- > 0 $, the norm
   \begin{equation}
      \left\Vert u\right\Vert =\inf\left\{  \lambda>0\left\vert \,\int_{\mathbb{R}^{N}}\left(  \left\vert \frac{\nabla u}{\lambda}\right\vert ^{z(x)}+M(x)\left\vert
      \frac{u}{\lambda}\right\vert ^{z(x)}\right) \leq1\right. \right\}  \label{n}
   \end{equation}
   is equivalent to norm $\left\Vert \, \cdot \,\right\Vert _{1,z(x)}$. If $ z_- > 1 $, the spaces $L^{z(x)}(\mathbb{R}^{N})$ and $W^{1,z(x)}(\mathbb{R}^{N})$ are reflexive
   and separable Banach spaces with these norms.

   \begin{proposition} \label{p1}
      The functional $\xi:W^{1,z(x)}(\mathbb{R}^{N}) \to \mathbb{R}$ defined by
   \begin{equation}
      \xi(u)=\int_{\mathbb{R}^N }\left(  \left\vert \nabla u\right\vert^{z(x)}+M(x)\left\vert u\right\vert ^{z(x)}\right) \text{,} \label{ps}
   \end{equation}
   has the following properties:

  \begin{enumerate}
    \item[\emph{(i)}] If $\left\Vert u\right\Vert \geq1$, then $\left\Vert u\right\Vert^{z_-}\leq\xi(u)\leq\left\Vert u\right\Vert ^{z_+}$.

    \item[\emph{(ii)}] If $\left\Vert u\right\Vert \leq1$, then $\left\Vert u\right\Vert^{z_+}\leq\xi(u)\leq\left\Vert u\right\Vert ^{z_-}$.
  \end{enumerate}
  In particular, for $ (u_n) \subset W^{ 1,z(x) }( \mathbb R^N ) $,
  \begin{gather*}
     \left\Vert u_n\right\Vert \rightarrow0 \iff \xi( u_n) \rightarrow0, \ \text{and}, \\
     ( u_n ) \ \text{is bounded in}\   W^{ 1,z(x) }( \mathbb R^N ) \iff \xi( u_n) \ \text{is bounded in} \ \mathbb R.
  \end{gather*}
\end{proposition}

\begin{remark} \label{r1}
   For the functional $\xi_z:L^{z(x)}(\mathbb{R}^{N})\rightarrow \mathbb{R}$ given by
   \[
      \xi_z(u)=\int_{\mathbb{R}^{N}}\left\vert u\right\vert^{z(x)}\text{,}
   \]
  the same conclusion of Proposition \ref{p1} also holds. Moreover, from $(i)$ and $(ii)$,
  \begin{equation}
     \left\vert u\right\vert _{z(x)}\leq \max\left\{\left(\int_{\mathbb{R}^{N}}\left\vert u\right\vert^{z(x)}\right) ^{1/z_{-}},\left(  \int_{\mathbb{R}^{N}}
     \left\vert u\right\vert ^{z(x)}\right)^{1/z_{+}} \right\} \text{.} \label{in}
\end{equation}

\end{remark}

Related to the Lebesgue space $L^{z(x)}(\mathbb{R}^{N})$, we have the following generalized H\"{o}lder's inequality.

\begin{proposition}
[{\cite[p.9]{Mu}}]\label{h} For $z\in L^{\infty}( \mathbb{R}^{N}) $ with $ z_->1$, let $z^{\prime }:\mathbb{R}^{N} \to \mathbb{R}$ be such that
\[
\frac{1}{z( x) }+\frac{1}{z^{\prime}( x) }=1 \ \text{a.e. in} \ \mathbb R^N.
\]
Then, for any $u\in L^{z( x) }( \mathbb{R}^{N})
$ and  $v\in L^{z^{\prime}( x) }( \mathbb{R}^{N}) $,
\begin{equation}
\left\vert \int_{\mathbb{R}^{N}}uv\,\right\vert
\leq \left( \frac{1}{z_-} + \frac{1}{z'_-} \right)\left\vert u\right\vert _{z( x) }\left\vert v\right\vert
_{z^{\prime}( x) }\text{.} \label{hi}
\end{equation}
\end{proposition}

\begin{proposition}
[{\cite[Theorems 1.1, 1.3]{FSZ}}]\label{pp1}Let $z:\mathbb{R}^{N} \to \mathbb{R}$ be a Lipschitz continuous satisfying
$ 1 < z_- \leq z_+ < N $ and $s:\mathbb{R}^{N} \rightarrow \mathbb{R}$ be a measurable function.

\begin{enumerate}
\item[\emph{(i)}] If $z\leq s\leq z^{\ast}$, the embedding $W^{1,z(x)}(\mathbb{R}^{N})\hookrightarrow L^{s(x)}(\mathbb{R}^{N})$ is continuous.

\item[\emph{(ii)}] If $z\leq s\ll z^{\ast}$, the embedding $W^{1,z(x)}(\mathbb{R}^{N})\hookrightarrow L_{\mathrm{loc}}^{s(x)}(\mathbb{R}^{N})$  is compact.
\end{enumerate}
\end{proposition}

The next two results are very important in our arguments and their proofs follows the same arguments explored in \cite{Kavian}, this form, we will omit their proofs.

\begin{proposition} [Brezis-Lieb's lemma, first version] \label{first Brezis-Lieb}
   Let $ z \in L^{\infty}(\mathbb R^N)$ with $ z_{-} \geq 1$ and $ ( \eta_n ) \subset L^{ z(x) } ( \mathbb R^N, \mathbb R^k ) $ verifying

   \begin{enumerate}
      \item [\emph{(i)}] $ \eta_n(x) \to \eta(x), \ \text{a.e. in} \ \mathbb R^N $;
      \item [\emph{(ii)}] $ \displaystyle \sup_{n \in \mathbb N } | \eta_n |_{ L^{ z(x) } ( \mathbb R^N, \mathbb R^k ) } < \infty $. \\
   \end{enumerate}
   Then, $ \eta \in L^{ z(x) } ( \mathbb R^N, \mathbb R^k ) $ and
   \begin{equation}
      \int_{ \mathbb R^N } \left( \left| \eta_n \right|^{ z(x) } - \left| \eta_n - \eta \right|^{ z(x) } - \left| \eta \right|^{ z(x) } \right) \, dx = o_n(1).
   \end{equation}
\end{proposition}

\begin{proposition} [Brezis-Lieb's lemma, second version] \label{second Brezis-Lieb}
   Let $ z \in L^{\infty}(\mathbb R^N)$ with $ z_{-} > 1$ and $ ( \eta_n ) \subset L^{ z(x) } ( \mathbb R^N, \mathbb R^k ) $ verifying

   \begin{enumerate}
      \item [\emph{(i)}] $ \eta_n(x) \to \eta(x), \ \text{a.e. in} \ \mathbb R^N $;
      \item [\emph{(ii)}] $ \displaystyle \sup_{n \in \mathbb N } | \eta_n |_{ L^{ z(x) } ( \mathbb R^N, \mathbb R^k ) } < \infty $. \\
   \end{enumerate}
   Then
   \begin{equation}
      \eta_n \rightharpoonup \eta \ \text{in} \ L^{ z(x) } ( \mathbb R^N, \mathbb R^k ).
   \end{equation}
\end{proposition}

The next proposition is a Brezis-Lieb type result and it applies an important role in our paper. For the case where $z$ is constant, the result is due to Alves \cite{Alves1} for $z \geq 2$ and Mercuri \& Willem \cite{MW} for $1<z<2$.

\begin{proposition} [Brezis-Lieb lemma, third version]
   Let  $ z \in L^{\infty}(\mathbb R^N)$ with $ z_{-} > 1 $ and $ ( \eta_n ) $ a sequence in $ L^{ z(x) }  ( \mathbb R^N, \mathbb R^k ) $ such that

   \begin{enumerate} 
      \item [\emph{(i)}] $ \eta_n(x) \to \eta(x), \ \text{a.e. in} \ \mathbb R^N $;
      \item [\emph{(ii)}] $ \displaystyle \sup_{n \in \mathbb N } | \eta_n |_{ L^{ z(x) } ( \mathbb R^N, \mathbb R^k ) } < \infty $. \\
   \end{enumerate}
   Then
   \begin{equation} \label{third Brezis-Lieb}
      \int_{\mathbb R^N} \left| \left| \eta_n \right|^{ z(x)-2 } \eta_n - \left| \eta_n - \eta \right|^{ z(x)-2 } \left( \eta_n - \eta \right) - \left| \eta \right|^{ z(x)-2 } \eta \right|^{ z'(x) }  = o_n(1).
   \end{equation}
 \end{proposition}

\noindent {\bf Proof.}
  In what follows, we set
   $$
      A(x,y) = | y |^{ z(x)-2 }y, \, \forall x \in \mathbb R^N, y \in \mathbb R^k.
   $$
  Our goal is to show that
   \begin{equation} \label{integral in region 1}
      \small{\int_{ \left\{ x \in \mathbb R^N; \, 1 < z(x) < 2 \right\} } \left| A \big( x, \eta_n(x) \big) - A \big( x, \eta_n(x) -  \eta(x) \big) - A \big( x, \eta(x) \big)  \right|^{ z'(x) } = o_n(1)}
   \end{equation}
   and
   \begin{equation} \label{integral in region 2}
       \small{\int_{ \left\{ x \in \mathbb R^N; \, z(x) \ge 2 \right\} } \left| A \big( x, \eta_n(x) \big) - A \big( x, \eta_n(x) -  \eta(x) \big) - A \big( x, \eta(x) \big)  \right|^{ z'(x) } = o_n(1),}
   \end{equation}
because if the above limits occur, we have that \eqref{third Brezis-Lieb} also occurs. This way, we will begin showing the limit \eqref{integral in region 1}. If the set    $ z^{ -1 } \big( (1,2) \big) $ has zero measure, we have nothing to do. Thereby, we will assume that  $ z^{ -1 } \big( (1,2) \big) $ has a positive measure and we will adapt the ideas found in \cite{MW}. First of all, we observe that
   \begin{equation} \label{sup in region 1}
      \alpha = \sup_{ { x \in z^{ -1 } ( (1,2) ) \atop y, h \in \mathbb R^k } \atop h \ne 0 } F(x,y,h) < \infty,
   \end{equation}
   where
   $$
      F(x,y,h) = \left| \frac{|y+h|^{ z(x)-2 }(y+h)-|y|^{ z(x)-2 }y}{|h|^{ z(x)-1 }} \right|.
   $$
   In fact, given any $ t > 0 $, it is easy to see that
   $$
      F(x,y,th) = F \left( x, \frac{y}{t}, h \right),
   $$
hence
   $$
      \alpha = \sup_{ { x \in z^{ -1 } ( (1,2) ) \atop y, h \in \mathbb R^k } \atop | h | = 1 } F(x,y,h).
   $$
   Firstly, if $ | y | \le 2 $, for any $ x \in z^{ -1 } \big( (1,2) \big) $, $ h \in \mathbb R^k $ with $ | h | = 1 $, it follows that
   $$
      \left| | y+h |^{ z(x)-2 }( y+h ) - | y |^{ z(x)-2 } y \right| \le 5,
   $$
 implying that
   \begin{equation} \label{sup in region 11}
      \alpha_1 = \sup_{ { x \in z^{ -1 } ( (1,2) ) \atop y, h \in \mathbb R^k } \atop | y | \le 2, | h | = 1 } F(x,y,h) < \infty.
   \end{equation}
   On the other hand, if $ | y | > 2 $, for any $ t \in [0,1] $ and $ h \in \mathbb R^k $ with $ | h | = 1 $, it holds
   $$
      | y+th | \ge | y | - t| h | > 1.
   $$
   Therefore, for each $ i = 1, \ldots, k $ and  $ x \in z^{ -1 } \big( (1,2) \big) $,
   \begin{align*}
      & \left| | y+h |^{ z(x)-2 }( y_i+h_i) - | y |^{ z(x)-2 }y_i \right| = \left| \int_0^1 \frac{d}{dt} | y+th |^{ z(x)-2 }(y_i+th_i) \, dt \right| \\
      & \phantom{ \left| \right. } = \left| \int_0^1 \left( | y + th |^{ z(x)-2 } h_i + \big( z(x)-2 \big) ( y_i + th_i ) | y + th |^{ z(x)-4 } ( y + th ) \cdot h \right) \, dt \right| \\
      & \phantom{ \left| \right. } \le \big( 3 - z(x) \big) \int_0^1 | y + th |^{ z(x)-2 } \, dt < 2 \int_0^1 1 \, dt = 2,
   \end{align*}
showing that
   \begin{equation} \label{sup in region 12}
      \alpha_2 = \sup_{ { x \in z^{ -1 } ( (1,2) ) \atop y, h \in \mathbb R^k } \atop | y | > 2, | h | = 1 } F(x,y,h) < \infty.
   \end{equation}
Combining \eqref{sup in region 11} with \eqref{sup in region 12}, we obtain \eqref{sup in region 1}.

A direct computation gives
   \begin{multline*}
      \left | A \big( x, \eta_n(x) \big) - A \big( x, \eta_n(x) - \eta(x) \big) - A \big( x, \eta(x) \big) \right| \\
      \le F \big( x, \eta_n(x) -  \eta(x), \eta(x) \big) | \eta(x) |^{ z(x)-1 } + | \eta(x) |^{ z(x)-1 } \le ( \alpha + 1 ) | \eta(x) |^{ z(x)-1 },
   \end{multline*}
   for all $ x \in z^{ -1 } \big( (1,2) \big) $, and so,
   $$
      \left | A \big( x, \eta_n(x) \big) - A \big( x, \eta_n(x) - \eta(x) \big) - A \big( x, \eta(x) \big) \right|^{ z'(x) } \le ( \alpha + 1 )^{ z'_+ } | \eta(x) |^{ z(x) },
   $$
   for all $ x \in z^{ -1 } \big( (1,2) \big) $, where $ z'(x) = \frac{z(x)}{z(x)-1}, \, \forall x \in \mathbb R^N $. Now, the limit \eqref{integral in region 1} follows from the last inequality together with Lebesgue's dominated convergence theorem.

In the proof of \eqref{integral in region 2}, we will adapt the ideas found in \cite{Alves1}.  If the set    $ z^{ -1 } \big( [2, \infty) \big) $ has zero measure, we have nothing to do. Thereby, we will assume that  $ z^{ -1 } \big( [2, \infty) \big) $ has a positive measure.  For each $ i = 1, \ldots, k $ and  $ x \in \mathbb R^N $, we have that
   \begin{multline*}
      A_i \big( x, \eta_n(x) \big) - A_i \big( x, \eta_n(x) - \eta(x) \big) \\
      = \left| \eta_n(x) \right|^{ z(x)-2 } \eta_n^i(x) - \left| \eta_n(x) -\eta(x) \right|^{ z(x)-2 } \left( \eta_n^i(x) - \eta_i(x) \right).
   \end{multline*}
   So, by the previous calculations,
   \begin{align*}
      { \textstyle \left| A_i \big( x, \eta_n(x) \big) \! - \! A_i \big( x, \eta_n(x) \! - \! \eta(x) \big) \right| } & { \textstyle \le \big( z(x) \! - \! 1 \big) | \eta(x) | { \displaystyle \int_0^1 } \negthickspace \left| \eta_n(x) \! + \! ( t \!- \! 1 ) \eta(x) \right|^{ z(x)-2 } dt } \\
      & \! \le  \big( z_+ - 1 \big) | \eta(x) | \left( \left| \eta_n(x) \right| + | \eta(x) |\right)^{ z(x)-2 }.
   \end{align*}
   Therefore
   $$
      \left| A \big( x, \eta_n(x) \big) - A \big( x, \eta_n(x) - \eta(x) \big) \right| \le C \left( | \eta(x) |^{ z(x)-1 } + | \eta(x) || \eta_n(x) |^{ z(x)-2 } \right),
   $$
   for all $ x \in z^{ -1 } \big( [2, \infty) \big) $. The above inequality combined with Young's inequality leads to
   $$
      \left| A \big( x, \eta_n(x) \big) - A \big( x, \eta_n(x) -  \eta(x) \big) \right| \le C( \epsilon ) | \eta(x) |^{ z(x)-1 } + \epsilon | \eta_n(x) |^{ z(x)-1 }, \forall \epsilon > 0.
   $$
Now, for each $ \epsilon > 0 $, $ n \in \mathbb N $, we define the function  $f_{ \epsilon, n }: \mathbb{R}^{N} \to \mathbb{R}$ given by
   $$
      f_{ \epsilon, n }(x) \! = \! \max \left\{ \left| A \big( x, \eta_n(x) \big) \! - \! A \big( x, \eta_n(x) \! - \! \eta(x) \big) \! - \! A \big( x, \eta(x) \big) \right| - \epsilon | \eta_n(x) |^{ z(x)-1 }, 0 \right\},
   $$
  which satisfies
$$
      f_{ \epsilon, n }(x) \to 0  \,\,\, \mbox{a.e in} \,\,\,  \ z^{ -1 } \big( [2, \infty) \big),  \ \text{as} \ n \to \infty,
$$
and
$$
      0 \le f_{ \epsilon, n }(x) \le \big( C ( \epsilon ) + 1 ) | \eta(x) |^{ z(x) -1 }, \, \forall x \in z^{ -1 } \big( [2, \infty) \big).
$$
    So, by Lebesgue's dominated convergence theorem,
    $$
       \int_{ z^{ -1 } \big( [2, \infty) \big) } f_{ \epsilon, n }^{ z'(x) } \to 0, \ \text{as} \ n \to \infty.
    $$
On the other hand, by the definition of $ f_{ \epsilon, n } $,
    $$
       \left| A \big( x, \eta_n(x) \big) - A \big( x, \eta_n(x) - \eta(x) \big) - A \big( x, \eta(x) \big) \right| \le \epsilon | \eta_n(x) |^{ z(x)-1 } + f_{ \epsilon, n}(x),
    $$
    for all $ x \in \mathbb R^N $. Consequently,
    \begin{multline*}
       \left| A \big( x, \eta_n(x) \big) - A \big( x, \eta_n(x) - \eta(x) \big) - A \big( x, \eta(x) \big) \right|^{ z'(x) } \\
       \le 2^{ z'_+} \left( \epsilon^{ z'_-} | \eta_n(x) |^{ z(x) } + f_{ \epsilon, n}^{ z'(x) } \right),
    \end{multline*}
    for all $ x \in \mathbb R^N $ and $ \epsilon > 0 $ sufficiently small. Thus,
    \begin{multline*}
       \varlimsup_n \int_{ z^{ -1 } \big( [2, \infty) \big) } \left| A \big( x, \eta_n(x) \big) - A \big( x, \eta_n(x) - \eta(x) \big) - A \big( x, \eta(x) \big) \right|^{ z'(x) } \\
       \le 2^{ z'_+} \epsilon^{ z'_-} \int_{ z^{ -1 } \big( [2, \infty) \big) }  | \eta_n(x) |^{ z(x) } \le C \epsilon^{ z'_-} , \, \forall \epsilon > 0 ,
    \end{multline*}
    which implies that \eqref{integral in region 2} holds. \qed

\section{Preliminary results}

In what follows, we will consider on $ W^{ 1,p(x) } \big( \mathbb{R}^{N} \big) $ the following norm
   $$
      \| u \| = \inf \left\{ \alpha > 0; \, \rho( \alpha^{ -1 }u ) \le 1 \right\},
   $$
   with
   $$
      \rho( u ) = \int_{\mathbb{R}^{N}}( | \nabla u |^{ p(x) } + V(x)| u |^{ p(x) }).
   $$
   Using well known arguments, we have that the energy functional \linebreak $I:W^{ 1,p(x) }( \mathbb{R}^{N} ) \to \mathbb{R} $ associated with $ (P) $, which is given by
  \begin{multline*}
     I( u ) = \int_{\mathbb{R}^{N}} \frac{1}{p(x)} \left( | \nabla u |^{ p(x) } + V(x) | u |^{ p(x) } \right) - \lambda \int_{\mathbb{R}^{N}} \frac{h(x)}{r(x) } ( u^+ )^{ r(x) } \\
      - \mu \int_{\mathbb{R}^{N}} \frac{1}{q(x)} ( u^+ )^{ q(x) } - \int_{\mathbb{R}^{N}} \frac{1}{p^{ \ast }(x) } ( u^+ )^{ p^{ \ast}(x) },
  \end{multline*}
 is well defined and $ I \in C^1 ( W^{ 1,p(x) }( \mathbb{R}^{N} ), \mathbb{R} ) $ with
   \begin{multline*}
      I'( u ) v = \int_{\mathbb{R}^{N}} (| \nabla u |^{ p(x)-2 } \nabla u \nabla v + V(x)| u |^{ p(x)-2 } u v ) -\lambda \int_{\mathbb{R}^{N}} h(x)( u^+ )^{ r(x)-1 }v \\
      - \mu \int_{\mathbb{R}^{N}} ( u^+ )^{ q(x)-1 }v - \int_{\mathbb{R}^{N}} ( u^+ )^{ p^{ \ast }(x)-1} v ,
   \end{multline*}
   for all $u, v \in W^{ 1,p(x) } ( \mathbb{R}^{N}). $

\begin{lemma}
   All $ (PS)_d $ sequences $ ( v_n ) $ for $ I $ are bounded. Furthermore, $ ( v_n^+ ) $ is a $ (PS)_d $ sequence for $ I $.
\end{lemma}

\noindent {\bf Proof.}   If there exist only a finite number of  terms $ (v_n) $ such that $ \rho( v_n ) > 1 $, then $ ( v_n ) $ is bounded and the proof is complete. Otherwise,
suppose the existence of a infinitely many terms of $ (v_n) $ such that $ \rho( v_n ) > 1 $.   Since $ ( v_n ) $ is a $ (PS)_d $ sequence, there is $n_0 \in \mathbb{N}$ such that
   $$
      I( v_n ) - \frac{1}{q_-} I'( v_n )v_n \le d+1 + \| v_n \|, \, n \geq n_0 .
   $$
On the other hand, using the fact that $ \rho( v_n ) > 1 $ and  H\"older's inequality, we get
   \begin{align*}
      I( v_n ) - \frac{1}{q_-} I'( v_n )v_n  & \ge \left( \frac{1}{p_+} - \frac{1}{q_-} \right) \| v_n \|^{ p_- } - \lambda \left( \frac{1}{r_-} - \frac{1}{q_-} \right) \int_{\mathbb{R}^{N}} h(x) | v_n |^{ r(x) } \\
      & \ge \left( \frac{1}{p_+} - \frac{1}{q_-} \right) \| v_n \|^{ p_- } - \lambda \left( \frac{1}{r_-} - \frac{1}{q_-} \right) C | h |_{ \Theta(x) } \big|  | v_n |^{ r(x) }  \big|_{ \frac{p^{\ast}(x)}{r(x)} },
   \end{align*}
and so,
   \begin{align*}
      I( v_n ) & - \frac{1}{q_-} I'( v_n )v_n \\
      & \ge \left( \frac{1}{p_+} - \frac{1}{q_-} \right) \| v_n \|^{ p_- } - \lambda \left( \frac{1}{r_-} - \frac{1}{q_-} \right) C | h |_{ \Theta(x) } \left( | v_n |_{ p^{ \ast }(x) }^{ r_- } + | v_n |_{ p^{ \ast }(x) }^{ r_+ } \right) \\
      & \ge \left( \frac{1}{p_+} - \frac{1}{q_-} \right) \| v_n \|^{ p_- } - \lambda \left( \frac{1}{r_-} - \frac{1}{q_-} \right) | h |_{ \Theta(x) } \left( C_1 \| v_n \|^{ r_- } + C_2 \| v_n \|^{ r_+ } \right).
   \end{align*}
From this, for $n \geq n_0$,
   \begin{multline*}
      d+1 + \|v_n \| \\
      \ge \left( \frac{1}{p_+} - \frac{1}{q_-} \right) \| v_n \|^{ p_- } - \lambda \left( \frac{1}{r_-} - \frac{1}{q_-} \right) | h |_{ \Theta(x) } \left( C_1 \| v_n \|^{ r_- } + C_2 \| v_n \|^{ r_+ } \right)
   \end{multline*}
which yields $ ( v_n ) $ is also bounded in this case.

   Now, we will prove that $(v_{n}^{+})$ is also a $(PS)_d$ sequence for $ I $. Note that the boundedness of $ ( v_n^- ) $ combined with the limit $\|I'(v_n)\| \to 0$ gives
$$
I'( v_n )v_n^- \to 0,
$$
from where it follows that
$$
\rho( v_n^- ) \to 0,
$$
or equivalently
$$
v_n^-  \to 0 \,\,\, \mbox{in} \,\,\, W^{1,p(x)}(\mathbb{R}^{N}).
$$
Now, a simple computation yields
$$
I( v_n)= I( v_n^+ ) + o_n(1) \,\,\, \mbox{and} \,\,\, I'(v_n)= I'( v_n^+ ) + o_n(1),
$$
proving that $ ( v_n^+ ) $ is a $ (PS)_d $ sequence.
\qed

\vspace{0.5 cm}

From the last lemma, hereafter we will assume that all $ (PS)_d $ sequences for $I$ are composed by nonnegative functions. Moreover, once that $W^{1,p(x)}(\mathbb{R}^N)$ is reflexive, if
$(v_n)$ is a $(PS)_d$ sequence for $ I $, we also assume that for some subsequence, still denoted by itself, there is $v \in W^{1,p(x)}(\mathbb{R}^{N})$ such that
$$
v_n \rightharpoonup v \,\,\, \mbox{in} \,\,\, W^{1,p(x)}(\mathbb{R}^{N}),
$$
$$
v_n(x) \to v(x) \,\,\, \mbox{a.e in } \,\,\, \mathbb{R}^{N},
$$
and
$$
v(x) \geq 0 \,\,\, \mbox{a.e in} \,\,\, \mathbb{R}^{N}.
$$
\vspace{0.5 cm}

The next lemma is a key point in our arguments, which can be found in \cite{AlvesFerreira}. However for the reader's convenience we will make its proof.

\vspace{0.5 cm}

\begin{lemma} \label{aew convergence}
   Let $ ( v_n ) $ be a $ (PS)_d $ sequence for $ I $ and $ v \in W^{ 1,p(x) } \big( \mathbb{R}^{N} \big) $ such that
   $ v_n \rightharpoonup v$ in $W^{ 1,p(x) } \big( \mathbb{R}^{N} \big) $. Then, $I'(v) = 0$. Hence, if $ v \ne 0 $, $ v $ is a nontrivial solution for $ (P) $.
\end{lemma}

\noindent {\bf Proof.}
Following a standard reasoning, it is sufficient to show that, up to a subsequence,
$$
\nabla v_n( x ) \to \nabla v( x ) \,\,\, \mbox{a.e in} \,\,\,  \ \mathbb{R}^{N}.
$$
We begin observing that, up to a subsequence, there exist two nonnegative measures  $ \mathfrak{m} $ and $ \mathfrak{n} $ in $ \mathcal{ M }\big( \mathbb{R}^{N} \big) $ such that
\begin{equation} \label{MED1}
| \nabla v_n |^{ p(x) }  \rightharpoonup \mathfrak{m} \ \text{in} \ \mathcal{M} \big( \mathbb{R}^{N} \big)
\end{equation}
and
\begin{equation} \label{MED2}
| v_n |^{ p^{ \ast }(x) }  \rightharpoonup \mathfrak{n} \ \text{in} \ \mathcal{M} \big( \mathbb{R}^{N} \big).
\end{equation}
In this case, according a concentration compactness principle in \cite{FuZa}, there exists a countable index set $ \mathfrak{I} $ such that
   \begin{gather*}
      \mathfrak{n} = | v |^{ p^{ \ast }(x) } \, dx + \sum_{ i \in \mathfrak{I} } \mathfrak{n}_i \delta_{ x_i }, \\
      \mathfrak{m} \ge | \nabla v |^{ p(x) } \, dx + \sum_{ i \in \mathfrak{I} } \mathfrak{m}_i \delta_{ x_i },
 \end{gather*}
and
$$
     \mathfrak{n}_i \le S \max \left\{ \mathfrak{m}_i^{ \frac{p^{ \ast}_+}{p_-} }, \mathfrak{m}_i^{ \frac{p^{ \ast}_-}{p_+} } \right\}
$$
where $ ( \mathfrak{n}_i )_{ i \in \mathfrak{I}}, ( \mathfrak{m}_i )_{ i \in \mathfrak{I} } \subset [0, \infty) $ and $ ( x_i )_{ i \in \mathfrak{I} } \subset \mathbb{R}^{N}. $ The constant $ S $ is given by
$$
      S = \sup_{ u \in W^{ 1,p(x) } ( \mathbb{R}^{N} ) \atop \| u \| \le 1 } \int_{\mathbb{R}^{N}} |u|^{ p^{ \ast }(x) }.
$$
Our first task is to prove that
$$
\mathfrak{m_i} = \mathfrak{n_i}, \, \forall i \in \mathfrak{I}.
$$
For this, let $ \varphi \in C^{ \infty }_0 \big( \mathbb{R}^{N} \big) $ such that
$$
\varphi(x) = 1 \,\,\,  \text{in} \ B_1( 0 ), \ \varphi(x) = 0  \ \text{in} \  B_2^c( 0 ) \ \text{and} \ 0 \le \varphi(x) \le 1 \, \forall \, x \in \mathbb{R}^{N}.
$$
Fixed $ i \in \mathfrak{I} $, we consider for each $ \epsilon > 0 $
$$
\varphi_{ \epsilon } (x) = \varphi \left( \frac{x - x_i}{\epsilon} \right) \,\,\, \forall x \in \mathbb{R}^{N}.
$$
Since $ ( v_n ) $ is bounded in $ W^{ 1,p(x) } \big( \mathbb{R}^{N} \big) $, the sequence $ ( \varphi_{ \epsilon } v_n ) $ is also bounded in $ W^{ 1,p(x) } \big( \mathbb{R}^{N} \big) $. Thus,
$$
      I'(v_n)( \varphi_{ \epsilon } v_n ) = o_n(1),
$$
that is,
   \begin{multline*}
      \int_{\mathbb{R}^{N}} (\varphi_{ \epsilon } | \nabla v_n |^{ p(x) } + v_n | \nabla v_n |^{ p(x)-2 } \nabla v_n  \nabla \varphi_{ \epsilon }) + \int_{\mathbb{R}^{N}} V(x) | v_n |^{ p(x) } \varphi_{ \epsilon }  \\
      = \lambda \int_{\mathbb{R}^{N}}h(x)|v_n|^{r(x)}\varphi_{ \epsilon }+ \mu \int_{\mathbb{R}^{N}} | v_n |^{ q(x) } \varphi_{ \epsilon }  + \int_{\mathbb{R}^{N}} | v_n |^{ p^{ \ast }(x) } \varphi_{ \epsilon } + o_n(1).
   \end{multline*}
   Taking the limits as $ n \to \infty $, the weak convergence  of $ (| \nabla v_n |^{ p(x) } ) $ and $ (| v_n |^{ p^{ \ast } (x) } )$ in $\mathcal {M}(\mathbb{R}^N)$ combined with
   the Lebesgue's dominated convergence theorem and Proposition \ref{second Brezis-Lieb}, give us
   \begin{multline} \label{epsilon limit}
      \int_{\mathbb{R}^{N}} \varphi_{ \epsilon } \, d \mathfrak{m} + \limsup_n \int_{\mathbb{R}^{N}} v_n | \nabla v_n |^{ p(x)-2 } \nabla v_n \nabla \varphi_{ \epsilon }  + \int_{\mathbb{R}^{N}} V(x) | v |^{ p(x) } \varphi_{ \epsilon }  \\
      = \lambda \int_{\mathbb{R}^{N}}h(x)|v|^{r(x)}\varphi_{ \epsilon }+  \mu \int_{\mathbb{R}^{N}} | v |^{ q(x) } \varphi_{ \epsilon }+ \int_{\mathbb{R}^{N}} \varphi_{ \epsilon } \, d \mathfrak{n}.
   \end{multline}
Using H\"older's inequality and the boundedness of $(v_n)$ in $W^{1,p(x)}(\mathbb{R}^{N})$,
   \begin{align*}
      & \left| \int_{\mathbb{R}^{N}} v_n | \nabla v_n |^{ p(x)-2 } \nabla v_n \cdot \nabla \varphi_{ \epsilon }  \right| \\
      & \phantom{ \left| \right. } \le \int_{\mathbb{R}^{N}} \left| \nabla v_n \right|^{ p(x)-1 } \left| v_n \nabla \varphi_{ \epsilon } \right|
        \le C \left| \left| \nabla v_n \right|^{ p(x)-1 } \right|_{ p'(x) } \big| v_n \left| \nabla \varphi_{ \epsilon } \right| \big|_{ p(x) } \\
      & \phantom{ \left| \right. } \le C \max \left\{ \left( \int_{\mathbb{R}^{N}} | v_n |^{ p(x) } \left| \nabla \varphi_{ \epsilon } \right|^{ p(x) }  \right)^{ \frac{1}{p_-} }, \left( \int_{\mathbb{R}^{N}} | v_n |^{ p(x) } \left| \nabla \varphi_{ \epsilon } \right|^{ p(x) } \right)^{ \frac{1}{p_+} } \right\},
   \end{align*}
where $ p'(x) = \frac{ p(x) }{ p(x)-1 }  \, \forall x \in \mathbb{R}^{N} $. Therefore, by Lebesgue's dominated convergence theorem,
   \begin{multline*}
      \limsup_n \left|\int_{\mathbb{R}^{N}} v_n | \nabla v_n |^{ p(x)-2 } \nabla v_n \cdot \nabla \varphi_{ \epsilon }  \right| \\
        \le C \max \left\{ \left( \int_{\mathbb{R}^{N}} | v |^{ p(x) } \left| \nabla \varphi_{ \epsilon } \right|^{ p(x) } \right)^{ \frac{1}{p_-} }, \left( \int_{\mathbb{R}^{N}} | v |^{ p(x) } \left| \nabla \varphi_{ \epsilon } \right|^{ p(x) }  \right)^{ \frac{1}{p_+} } \right\}.
   \end{multline*}
Furthermore, by H\"older's inequality
   $$
      \int_{\mathbb{R}^{N}} | v |^{ p(x) } \left| \nabla \varphi_{ \epsilon } \right|^{ p(x) }
        \le C \left| | v |^{ p(x) } \right|_{ L^{ \frac{N}{N-p(x)} } \big( B_{ 2 \epsilon }(x_i) \big) } \big| \left| \nabla \varphi_{ \epsilon } \right|^{ p(x) } \big|_{ L^{ \frac{N}{p(x)} } \big( B_{ 2 \epsilon }(x_i) \big) }. \\
   $$
   Once that
   $$
      \int_{ B_{ 2 \epsilon }(x_i) } \left| \nabla \varphi_{ \epsilon } \right|^N  = \int_{ B_2(0) } \left| \nabla \varphi \right|^N ,
   $$
   we derive
   \begin{multline*}
      \big| \left| \nabla \varphi_{ \epsilon } \right|^{ p(x) } \big|_{ L^{ \frac{N}{p(x)} } \big( B_{ 2 \epsilon }(x_i) \big) } \\
      \phantom{ \big| } \le \max \left\{ \left( \int_{ B_{ 2 \epsilon }(x_i) }  \left| \nabla \varphi_{ \epsilon } \right|^N  \right)^{ \frac{1}{\left( \frac{N}{p} \right)_-}  }, \left( \int_{ B_{ 2 \epsilon }(x_i) }  \left| \nabla \varphi_{ \epsilon } \right|^N  \right)^{ \frac{1}{\left( \frac{N}{p} \right)_+}  } \right\} \leq C
   \end{multline*}
 for some positive constant $C$, which is independent of $\epsilon$.  Thereby,
   $$
      \int_{\mathbb{R}^{N}} | v |^{ p(x) } \left| \nabla \varphi_{ \epsilon } \right|^{ p(x) } \le C \left| | v |^{ p(x) } \right|_{ L^{ \frac{N}{N-p(x)} } \big( B_{ 2 \epsilon }(x_i) \big) },
   $$
   and so
   \begin{multline*}
     \limsup_n \left|\int_{\mathbb{R}^{N}} v_n | \nabla v_n |^{ p(x)-2 } \nabla v_n \cdot \nabla \varphi_{ \epsilon } \right| \\
      \le C \max \left\{ \left| | v |^{ p(x) } \right|_{ L^{ \frac{N}{N-p(x)} } \big( B_{ 2 \epsilon }(x_i) \big) }^{ \frac{1}{p_-} }, \left| | v |^{ p(x) } \right|_{ L^{ \frac{N}{N-p(x)} } \big( B_{ 2 \epsilon }(x_i) \big) }^{ \frac{1}{p_+} } \right\}. \\
   \end{multline*}
   But,
   \begin{multline*}
      \left| | v |^{ p(x) } \right|_{ L^{ \frac{N}{N-p(x)} } \big( B_{ 2 \epsilon }(x_i) \big) } \\
      \le \max \left\{ \left( \int_{ B_{ 2 \epsilon }(x_i) }  | v |^{ p^{ \ast }(x) } \right)^{ \frac{1}{\left( \frac{N}{N-p} \right)_- } }, \left( \int_{ B_{ 2 \epsilon }(x_i) }  | v |^{ p^{ \ast }(x) }  \right)^{ \frac{1}{\left( \frac{N}{N-p} \right)_+ } } \right\}
   \end{multline*}
from where it follows that
   $$
      \lim_{\epsilon \to 0}\limsup_n \left| \int_{\mathbb{R}^{N}} v_n | \nabla v_n |^{ p(x)-2 } \nabla v_n \nabla \varphi_{ \epsilon } \right| = 0
   $$
implying that
$$
      \lim_{\epsilon \to 0}\limsup_n  \int_{\mathbb{R}^{N}} v_n | \nabla v_n |^{ p(x)-2 } \nabla v_n \nabla \varphi_{ \epsilon } = 0.
 $$
Now, taking the limit as $ \epsilon \to 0 $ in  \eqref{epsilon limit},  we get
   \begin{equation} \label{m_i = n_i}
      \mathfrak{m}_i = \mathfrak{m}( x_i ) = \mathfrak{n}( x_i ) = \mathfrak{n}_i.
   \end{equation}
Once that
$$
\frac{ p^{\ast}_- }{ p_+ } \le \frac{ p^{ \ast }_+ }{ p_- },
$$
we have that
\begin{equation} \label{eta1}
      \mathfrak{n}_i^{ \frac{p_+}{p^{ \ast }_-} } \le \left( S^{ \frac{p_+}{p^{ \ast }_-} } + S^{ \frac{p_-}{{p^{ \ast }_+}} } \right) \mathfrak{m}_i, \ \text{if} \,\,\, \mathfrak{m}_i <1
 \end{equation}
and
\begin{equation} \label{eta2}
\mathfrak{n}_i^{ \frac{p_-}{p^{ \ast }_+} } \le \left( S^{ \frac{p_+}{p^{ \ast }_-}} + S^{ \frac{p_-}{{p^{ \ast }_+}} } \right) \mathfrak{m}_i \,\,\, \mbox{if} \,\,\, \mathfrak{m}_i \geq 1.
\end{equation}
Thus, from \eqref{m_i = n_i} - \eqref{eta2}, if $ \mathfrak{n}_i > 0 $ for some $ i \in \mathfrak{I} $, there exists $\alpha>0$, which is independent of $ i $, such that
\begin{equation} \label{eta3}
\mathfrak{n}_i \geq \alpha.
\end{equation}
Recalling that
  \begin{equation} \label{n_i sum}
     \sum_{ i \in \mathfrak{I} \atop \mathfrak{m}_i < 1 } \mathfrak{n}_i^{ \frac{ p_+ }{ p^{\ast}_- } } +
       \sum_{ i \in \mathfrak{I} \atop \mathfrak{m}_i \ge 1 } \mathfrak{n}_i^{ \frac{ p_- }{ p^{\ast}_+ } } \le C \sum_{ i \in \mathfrak{I} } \mathfrak{m}_i < \infty,
   \end{equation}
the inequality \eqref{eta3} gives $ \tilde{\mathfrak{I}}= \left\{ i \in \mathfrak{I}; \, \mathfrak{n}_i > 0 \right\}$ is a finite set.  From this, one of the two possibilities below occurs: \\

  \noindent $a)$ There exist $ \mathfrak{n}_{i_{1}}, \ldots, {\mathfrak n}_{i_{s}} > 0 $ for a maximal $ s \in \mathbb{N} $; \\
  \mbox{}\\
  \noindent $b)$ $ \mathfrak{n}_i = 0 $, for all $ i \in \mathfrak{I}  $. \\

   We begin analyzing $a)$. For this, choose $ 0 < \epsilon_0 < 1 $ sufficiently small such that
   $$
      B_{ \epsilon_0 }(x_1), \cdots, B_{ \epsilon_0 }(x_s) \subset B_{ \frac{1}{\epsilon_0} }(0) \ \text{and} \ B_{ \epsilon_0 }(x_i) \cap B_{ \epsilon_0 }(x_j) = \emptyset, \ i \ne j,
   $$
  where $ x_1, \ldots, x_s $ are the singular points related to $ \mathfrak{n}_{i_1}, \ldots, \mathfrak{n}_{i_s}$, respectively. We set
   $$
      \psi_{ \epsilon }(x) = \varphi( \epsilon x ) - \sum_{ i=1 }^s \varphi \left( \frac{x-x_i}{\epsilon} \right) \, \forall x \in \mathbb{R}^{N}.
   $$
   Then, for $ 0 < \epsilon < \frac{1}{2} \epsilon_0 $,
   $$
      \psi_{ \epsilon}(x) =
      \begin{cases}
         0, & \ \text{if} \ x \in \displaystyle \bigcup_{ i=1 }^s B_{ \frac{\epsilon}{2} }( x_i ) \\
         1, & \ \text{if} \ x \in A_{ \epsilon } = B_{ \frac{1}{\epsilon} }( 0 ) \setminus \displaystyle \bigcup_{ i=1 }^s B_{ 2 \epsilon }( x_i ),
      \end{cases}
   $$
and 
   $$
      \text{supp} \, \psi_{ \epsilon } \subset \overline{B_{ \frac{2}{\epsilon} }( 0 )} \setminus \bigcup_{ i=1 }^s B_{ \frac{\epsilon}{2} }( x_i )
   $$
loading to 
   $$
      \int_{\mathbb{R}^{N}} | v_n |^{ p^{ \ast }(x) } \psi_{ \epsilon }  \to \int_{\mathbb{R}^{N}} | v |^{ p^{ \ast }(x) } \psi_{ \epsilon }.
   $$
   Since
   $$
      I'( v_n )( v_n \psi_{ \epsilon } ) = o_n(1) \ \text{and} \ I'( v_n )( v \psi_{ \epsilon } ) = o_n(1),
   $$
   repeating the same type of arguments for the case where the exponents are constant, we obtain
   $$
      \lim_n \int_{ A_{ \epsilon } } (P_n + V(x)Q_n)=0,
   $$
   where
 $$
 P_n( x ) = \left( \left| \nabla v_n \right|^{ p(x)-2 } \nabla v_n - \left| \nabla v \right|^{ p(x)-2 } \nabla v \right)\Big( \nabla v_n - \nabla v \Big) \,\,\, \forall x \in \mathbb{R}^{N} \,\,\, \mbox{and} \,\,\, \forall n \in \mathbb{N}.
  $$
and
$$
Q_n(x)=\left( \left| v_n \right|^{ p(x)-2 } v_n - \left| v \right|^{ p(x)-2 } v \right)\Big( v_n - v \Big) \,\,\, \forall x \in \mathbb{R}^{N} \,\,\, \mbox{and} \,\,\, \forall n \in \mathbb{N}.
$$
Since
   \begin{equation}
      P_n(x) \ge \label{P_n inequality}
      \begin{cases}
         \frac{2^{ 3-p_+ }}{p_+} \left| \nabla v_n - \nabla v \right|^{ p(x) }, \ \text{if} \ p(x) \ge 2 \\
         \left( p_- - 1 \right) \frac{\left| \nabla v_n - \nabla v \right|^2}{{\left( \left| \nabla v_n \right| + \left| \nabla v \right| \right)}^{ 2-p(x) }}, \ \text{if} \ 1< p(x) < 2,
      \end{cases}
   \end{equation}
it follows that
   $$
      \int_{ A_{ \epsilon } } P_n \ge C \int_{ A_{ \epsilon } \cap \left\{ x \in \mathbb{R}^{N}; \, p(x) \ge 2 \right\} } \left| \nabla v_n - \nabla v \right|^{ p(x) } \ge 0.
   $$
 Thus,
  \begin{equation} \label{P1}
      \lim_n \int_{  A_{ \epsilon } \cap \left\{ x \in \mathbb{R}^{N}; \, p(x) \ge 2 \right\} } \left| \nabla v_n - \nabla v \right|^{ p(x) } = 0.
  \end{equation}
On the other hand, by H\"older's inequality
   \begin{align*}
      & \int_{  A_{ \epsilon } \cap \left\{ x \in \mathbb{R}^{N}; \, 1 <  p(x) < 2 \right\} } \left| \nabla v_n - \nabla v \right|^{ p(x) } \\
      & \phantom{ \int_{\mathbb{R}^{N}} } \le C \left| \frac{\left| \nabla v_n - \nabla v \right|^{ p(x) }}{\big( \left| \nabla v_n \right| + \left| \nabla v \right| \big)^{ \frac{p(x)(2-p(x))}{2} } } \right|_{ L^{ \frac{2}{p(x)} }  \left( \tilde{A_{ \epsilon }} \right)  }
         \! \! \! \left| \big( \left| \nabla v_n \right| + \left| \nabla v \right| \big)^{ \frac{p(x)(2-p(x))}{2} } \right|_{ L^{ \frac{2}{2-p(x)} } \left( \tilde{A_{ \epsilon }} \right) },
  \end{align*}
   where $ \tilde{A_{ \epsilon }} = A_{ \epsilon } \cap \left\{ x \in \mathbb{R}^{N}; \, 1 <  p(x) < 2 \right\} $. From relation \eqref{P_n inequality}, the right side of above inequality goes to zero. Hence,
\begin{equation} \label{P2}
\lim_{n}\int_{  A_{ \epsilon } \cap \left\{ x \in \mathbb{R}^{N}; \, 1 <  p(x) < 2 \right\} } \left| \nabla v_n - \nabla v \right|^{ p(x) }=0.   \\
\end{equation}
Now  (\ref{P1}) combined with (\ref{P2}) gives
$$
\lim_{n}\int_{  A_{ \epsilon } } \left| \nabla v_n - \nabla v \right|^{ p(x) }=0.
$$
The same arguments can be used to prove that
$$
\lim_{n}\int_{  A_{ \epsilon } } V(x)\left| v_n - v \right|^{ p(x) }=0.
$$
Therefore,
$$
v_n \to v \,\,\, \mbox{in} \,\,\, W^{1,p(x)}(A_{ \epsilon }).
$$
The last limit yields, up to a subsequence,
   $$
       \nabla v_n( x ) \to \nabla v( x ) \ \,\,\, \mbox{a.e in} \,\,\,  \ A_{\epsilon} \,\,\, ( 0 < \epsilon < \frac{1}{2} \epsilon_0).
   $$
Observing that
   $$
      \mathbb{R}^{N} \setminus \left\{ x_1, x_2, \ldots, x_s \right\} = \bigcup_{ n \in \mathbb{N} \atop \frac{1}{n} < \frac{1}{2} \epsilon_0 } A_{ \frac{1}{n} },
   $$
we conclude  by a diagonal argument, that there is a subsequence of $(v_n)$, still denoted by itself, such that
   $$
      \nabla v_n( x ) \to \nabla v( x ) \,\,\,  \mbox{a.e in} \,\, \, \mathbb{R}^{N}.
   $$
For the case $b)$, we consider
   $$
      \psi_{ \epsilon }( x ) = \varphi( \epsilon x ) \,\,\, \forall x \in \mathbb{R}^{N} \,\,\,  \text{and} \,\,  A_{ \epsilon } = B_{ \frac{1}{\epsilon} }(0), \ \epsilon > 0.
   $$
Repeating the same arguments used in the case $a)$, we have that
$$
v_n \to v \,\,\, \mbox{in} \,\,\, W^{1,p(x)}(B_{ \frac{1}{\epsilon} }(0)).
$$
This way, there is again a subsequence of $(v_n)$, still denoted by itself, such that
$$
      \nabla v_n( x ) \to \nabla v( x ) \ \,\,\,  \mbox{a.e in} \,\, \, \ \mathbb{R}^{N}.
$$
\qed

\begin{lemma} \label{estimate}
Let $ ( v_n ) $  be a $ (PS)_d $ sequence  for $ I $ with  $ v_n \rightharpoonup v $ in $ W^{ 1,p(x) } \big( \mathbb{R}^{N} \big) $. Then, there exists a constant
   $ M > 0 $, which is independent of $ \lambda $ and $ \mu $,  such that
   $$
      I(v) \ge -M \left( \lambda^{ \Theta_- } + \lambda^{ \Theta_+ } \right).
   $$
\end{lemma}

\noindent {\bf Proof.} \, From Lemma \ref{aew convergence}, $I'(v)v=0$, or equivalently ,
$$
       \int_{\mathbb{R}^{N}} | \nabla v|^{ p(x) } + V(x) v^{ p(x) } =  \lambda \int_{\mathbb{R}^{N}} h(x) v^{ r(x) } + \mu \int_{\mathbb{R}^{N}} v^{ q(x) }  + \int_{\mathbb{R}^{N}} v^{ p^{ \ast }(x) } .
$$
From this,
$$
I(v)  \ge \lambda \left( \frac{1}{p_+} - \frac{1}{r_-} \right) \int_{\mathbb{R}^{N}} h(x) v^{ r(x) }+ \left( \frac{1}{p_+} - \frac{1}{p^{ \ast }_-} \right) \int_{\mathbb{R}^{N}} v^{ p^{ \ast }(x) },
$$
which together with Young's inequality implies that for all  $ \epsilon > 0 $,
   \begin{align*}
      I(v) & \ge \epsilon \left( \frac{1}{p_+} - \frac{1}{r_-} \right) \int_{\mathbb{R}^{N}} v^{ p^{ \ast }(x) }+ C( \epsilon, x) \left( \frac{1}{p_+} - \frac{1}{r_-} \right) \int_{\mathbb{R}^{N}} \lambda^{ \Theta(x) } h^{ \Theta(x) } \\
      & \phantom{ \ge } + \left( \frac{1}{p_+} - \frac{1}{p^{ \ast }_-} \right) \int_{\mathbb{R}^{N}} v^{ p^{ \ast }(x) },
      \end{align*}
  where
   $$
      C( \epsilon, x ) = \frac{1}{\Theta(x) \left( \epsilon \frac{p^{ \ast }(x)}{r(x)} \right)^{ \frac{r(x)\Theta(x)}{p^{ \ast }(x)} } }.
   $$
Fixing
   $$
     0< \epsilon < \min \left\{1, \left( \frac{1}{r_-}  - \frac{1}{p_+} \right)^{-1}\left( \frac{1}{p_+} - \frac{1}{p^{ \ast }_-} \right) \right\},
   $$
it follows that
   $$
      I(u) \ge -M \left( \lambda^{ \Theta_- } + \lambda^{ \Theta_+ } \right),
   $$
   where
   $$
     M = \frac{1}{\Theta_- \epsilon^{ \Theta_+ -1 } } \left( \frac{1}{r_-} - \frac{1}{p_+} \right) \int_{\mathbb{R}^{N}} h^{ \Theta(x) }.
   $$
\qed

The next result is an important step to understand the behavior of the $(PS)$ sequences of  $I$.

\begin{lemma} \label{decomposicao}
   Let $ ( v_n ) $ be  a bounded sequence in $ W^{ 1,p(x) }( \mathbb R^N ) $ such that $ v_n(x) \to v(x) \ \text{and} \ \nabla v_n(x) \to \nabla v(x) \,\,\, \mbox{a.e in} \,\,\, \mathbb R^N $. Then,
   $$
      I(v_n)-I_\mu(v_n-v)-I(v)=o_n(1) \leqno{i)}
   $$
   and
   $$
      I'(v_n)-I_\mu '(v_n-v)-I'(v)=o_n(1). \leqno{ii)}
   $$
   Consequently, if $ (v_n) $ is a $ (PS)_d $ sequence for $ I $ with weak limit $ v \in W^{ 1,p(x) }( \mathbb R^N ) $, setting $ w_n = v_n - v $, we have that for some subsequence, $ (w_n) $ is a $ (PS)_{d-I(v)} $ sequence for $ I_\mu$.
\end{lemma}

\noindent {\bf Proof.}   From definitions of  $I$ and $I_\mu$, we derive that
      \begin{align*}
      I & ( v_n ) - I_\mu( v_n-v ) - I(v)= \\
      & = \int_{\mathbb{R}^{N}} \frac{1}{p(x)} \left( \left| \nabla v_n \right|^{ p(x) } - \left| \nabla v_n - \nabla u \right|^{ p(x) } - \left| \nabla v \right|^{ p(x) } \right)  \\
      & { \textstyle \phantom{ = } + { \displaystyle \int_{\mathbb{R}^{N}} } \negthickspace \frac{V(x)}{p(x)} \negthickspace \left( \! v_n^{ p(x) } \negthickspace - \left| v_n - v \right|^{ p(x) } \negthickspace - v^{ p(x) }  \! \right)
           - \mu { \displaystyle \int_{\mathbb{R}^{N}} } \negthickspace \frac{1}{q(x)} \negthickspace \left( \! v_n^{ q(x) } \negthickspace - \left| v_n - v \right|^{ q(x) } \negthickspace - v^{ q(x) } \! \right) } \\
      & \phantom{ = } - \int_{\mathbb{R}^{N}} \frac{1}{ p^{ \ast }(x) } \left( v_n^{ p^{ \ast }(x) } - \left| v_n - v \right|^{ p^{ \ast }(x) } - v^{ p^{ \ast }(x) } \right) 
           - \lambda \int_{\mathbb{R}^{N}} \frac{h(x)}{r(x)} \left( v_n^{ r(x) } - v^{ r(x) } \right) .
   \end{align*}
By Propositions \ref{first Brezis-Lieb} and \ref{second Brezis-Lieb}, we observe that the right side of the last inequality is $ o_n(1) $, and so,
   $$
      I( v_n ) - I_\mu( v_n-v ) - I( v ) = o_n(1),
   $$
 showing $i)$.

   Now, to prove $ii)$, we fix $ \varphi \in W^{ 1, p(x) } \big( \mathbb{R}^{N} \big) $ with $ \| \varphi \| =1 $. Using H\"older's inequality together with Sobolev's embedding, it follows that there is a positive  constant $C$ such that
  $$
  \left[ I'( v_n ) - I'_\mu( v_n-v ) - I'( v ) \right] \varphi \big|  \leq C( A_1(n)+A_2(n)+A_3(n)+A_4(n)+A_5(n))
  $$
 where
$$
 \begin{array}{l}
 A_1(n)=\left| \left| \nabla v_n \right|^{ p(x)-2 } \nabla v_n - \left| \nabla v_n - \nabla v \right|^{ p(x)-2 } \left( \nabla v_n - \nabla v \right) - \left| \nabla v \right|^{ p(x)-2 } \nabla v \right|_{ p'(x) }, \\
 A_2(n)=\left| \ v_n^{ p(x)-2 } v_n \, - \, \left| v_n-v \right|^{ p(x)-2 } \left( v_n-v \right) \, - \, v^{ p(x)-2 } v \ \right|_{ p'(x)}, \\
 A_3(n)= \mu \left| \ v_n^{ q(x)-2 } v_n - \left| v_n-v \right|^{ q(x)-2 } \left( v_n-v \right) - v^{ q(x)-2 } v \ \right|_{ q'(x) } ,\\
 A_4(n)= \left| v_n^{ p^{ \ast }(x)-2 } v_n - \left| v_n-v \right|^{ p^{ \ast }(x)-2 } \left( v_n-v \right) - v^{ p^{ \ast }(x)-2 } v \right|_{ {p^{ \ast }}'(x) },
 \end{array}
$$
and
$$
\begin{array}{l}
A_5(n)= \lambda \displaystyle \int_{\mathbb{R}^{N}} h(x) \left| \left( v_n^{ r(x)-1 } - v^{ r(x)-1 } \right) \varphi \right|.
\end{array}
 $$
From Proposition \ref{third Brezis-Lieb}, $A_i(n)=o_n(1)$ for $i=1,2,3,4$. Related to $A_5(n)$, we have that
   $$
      \int_{\mathbb{R}^{N}} h(x) \left| \left( v_n^{ r(x)-1 }  - v^{ r(x)-1 } \right) \varphi \right|= \int_{\mathbb{R}^{N}} h^{ \frac{1}{r'(x)} } \left| v_n^{ r(x)-1 } - v^{ r(x)-1 } \right| h^{ \frac{1}{r(x)} } | \varphi |.
   $$
   Since
   $$
      h^{ \frac{1}{r'(x)} } \left| v_n^{ r(x)-1 } - v^{ r(x)-1 } \right| \in L^{ r'(x) } \big( \mathbb{R}^{N} \big) \ \text{and} \ h^{ \frac{1}{r(x)} } | \varphi | \in L^{ r(x) } \big( \mathbb{R}^{N} \big),
   $$
   by H\"older's inequality,
   $$
      \int_{\mathbb{R}^{N}} h(x) \left| \left( v_n^{ r(x)-1 }  - v^{ r(x)-1 } \right) \varphi \right| \le C \bigg| h^{ \frac{1}{r'(x)} } \left| v_n^{ r(x)-1 } - v^{ r(x)-1 } \right| \bigg|_{ r'(x) }.
   $$
Now, our goal is to prove that
   $$
      \bigg| h^{ \frac{1}{r'(x)} } \left| v_n^{ r(x)-1 } - v^{ r(x)-1 } \right| \bigg|_{ r'(x) } \to 0,
   $$
  or equivalently,
\begin{equation} \label{EQ1}
     \int_{\mathbb{R}^{N}} h(x) \left| v_n^{ r(x)-1 } - v^{ r(x)-1 } \right|^{ r'(x) } \to 0.
\end{equation}
To this end, we define
   $$
      V_n( x ) = \left| v_n^{ r(x)-1 } - v^{ r(x)-1 } \right|^{ r'(x) } \,\,\,  \forall n \in \mathbb{N}.
   $$
   Then,
   $$
      V_n( x ) \to 0 \ \,\,\,  \mbox{a.e in} \,\, \, \ \mathbb{R}^{N},
   $$
   and $( V_n )$ is bounded in $ L^{ \frac{p^{ \ast }(x)}{r(x)} } \big( \mathbb{R}^{N} \big)$. Therefore, by Proposition \ref{second Brezis-Lieb}, it follows that
  $$
  V_n \rightharpoonup 0 \,\,\, \mbox{in} \,\,\, L^{ \frac{p^{ \ast }(x)}{r(x)} } \big( \mathbb{R}^{N} \big).
  $$
Thus,
$$
\int_{\mathbb{R}^{N}}h(x)V_n(x) \to 0
$$
proving \eqref{EQ1}. Consequently,
   $$
      \big\| I'( v_n ) - I'_\mu( v_n-v ) - I'( v ) \big\| = o_n(1),
   $$
or yet
   $$
      I'( v_n ) - I'_\mu( v_n-v ) - I'( v ) = o_n(1),
   $$
finishing the proof. \qed

\begin{lemma} \label{compactness}
   Suppose $ \mu \ge \mu_0 $, where $ \mu_0 $ is given in \ref{ESTIMATIVAPASSO}). Then, I verifies the $ (PS)_d $ condition for
   $$
      d < c_\mu- M \left( \lambda^{ \Theta_- } + \lambda^{ \Theta_+ } \right).
   $$
\end{lemma}

\noindent {\bf Proof.} \,  Let $ ( v_n ) $ be a $ (PS)_d $ sequence for $ I $ with $ d $ as above. We know that there exists $ v \in W^{ 1,p(x) } \big( \mathbb{R}^{N} \big) $ such that
$$
      v_n \rightharpoonup v \ \text{in} \ W^{ 1,p(x) } \big( \mathbb{R}^{N} \big) ,
$$
and
$$
      v_n(x) \to v(x), \ \,\,\,  \mbox{a.e in} \,\, \, \ \mathbb{R}^{N}.
$$
 Setting $ w_n = v_n - v $, by Lemma \ref{decomposicao}, we see that $ ( w_n ) $ is a $ (PS)_{ d-I(v) } $ sequence for $ I_\mu $.  Thus, up to a subsequence, we can assume that
   $$
      \int_{\mathbb{R}^{N}}( | \nabla w_n |^{ p(x) } + V(x) | w_n |^{ p(x) }) \to L \ge 0.
   $$
Next, we will show that $ L = 0 $. To this end, we recall that only one of the below possibilities hold: \\

\noindent a) \, There is $ R > 0 $ such that
   $$
      \lim_{n} \sup_{ y \in \mathbb{R}^{N} } \int_{ B_R(y) } | w_n |^{ p(x) } = 0
   $$
or \\

\noindent b) \, For each $R>0$, there are $ \eta >0$, a subsequence of $(w_n)$, still denoted by itself,   and $ ( y_n ) \subset \mathbb{R}^{N} $ \big( which we can suppose in $ \mathbb{Z}^{N} $ \big) such that
   $$
      \varlimsup_n \int_{ B_R( y_n ) } | w_n |^{ p(x) } \ge \eta.
   $$
We will show that $b)$ does not hold. Arguing by contradiction, if $b)$ is true, we  define
   $$
     \widehat{ w_n }(x) = w_n( x+y_n ), \ x \in \mathbb{R}^{N}.
   $$
   Then, by a simple computation,
   $$
      \ I_\mu( \widehat{ w_n } ) = I_\mu( w_n ) \ \text{and} \ I'_\mu( \widehat{ w_n } ) = o_n(1).
   $$
   So, $ ( \widehat{ w_n } ) $ is also a $ (PS)_{ d - I(v) } $ sequence for $ I_\mu$. Let $ \widehat{ w } \in W^{ 1,p(x) }( \mathbb{R}^{N} ) \setminus \left\{ 0 \right\} $ the weak limit of $ \widehat{ w_n } $.
   Since $ I'_\mu( \widehat{ w } ) = 0 $ and $\widehat{ w }  \not=0$,  it follows from the definition of $c_\mu$ that
   \begin{align*}
      c_\mu &  \le I_\mu( \widehat{ w} ) = I_\mu( \widehat{w } ) - \frac{1}{p_+} I'_\mu( \widehat{w } )\widehat{ w } \\
      & \le \varliminf_{ n } \left( \int_{\mathbb{R}^{N}} \left( \frac{1}{p(x)} - \frac{1}{p_+} \right) \left( | \nabla \widehat{ w_n } |^{ p(x) } + V(x) \left| \widehat{ w_n } \right|^{ p(x) } \right) \right.\\
      & \left. \phantom{ \le } + \mu \int_{\mathbb{R}^{N}} \left( \frac{1}{p_+} - \frac{1}{q(x)} \right) \left| \widehat{ w_n } \right|^{ q(x) } + \int_{\mathbb{R}^{N}} \left( \frac{1}{p_+} - \frac{1}{p^{ \ast }(x)} \right) \left| \widehat{ w_n } \right|^{ p^{ \ast }(x) } \right) \\
      & = \varliminf_{ n } \left( I_\mu( \widehat{ w_n } ) - \frac{1}{p_+} I'_\mu( \widehat{ w_n } )\widehat{ w_n } \right) = d - I(v) \le d + M \left( \lambda^{ \Theta_- } + \lambda^{ \Theta_+ } \right).
   \end{align*}
   Thus,
   $$
      c_\mu - M \left( \lambda^{ \Theta_- } + \lambda^{ \Theta_+ } \right) \le d,
   $$
   which is a contradiction with the hypothesis on $ d $. Therefore, $ b) $ does not hold. Then $ a) $ holds,  and  by Lemma $ 3.1 $ in \cite{FZZ},
   $$
      w_n \to 0 \ \text{in} \ L^{ q(x) } \big( \mathbb{R}^{N} \big),
   $$
 or equivalently,
\begin{equation}\label{L1}
      \int_{\mathbb{R}^{N}} | w_n |^{ q(x) }  \to 0.
\end{equation}
 Since $ I'_\mu( w_n )w_n = o_n(1) $, we derive that
   $$
      \int_{\mathbb{R}^{N}} | w_n |^{ p^{ \ast }(x) }  \to L.
   $$
By (\ref{L1}),
   \begin{align*}
       d - I(v) + o_n(1) & = I_\mu( w_n ) + \mu \int_{\mathbb{R}^{N}} \frac{1}{q(x)} | w_n |^{ q(x) }  \\
      & = \int_{\mathbb{R}^{N}} \frac{1}{p(x)} \left( | \nabla w_n |^{ p(x) } + V(x) | w_n |^{ p(x) } \right)  - \int_{\mathbb{R}^{N}} \frac{1}{p^{ \ast }(x)} | w_n |^{ p^{ \ast }(x) }
   \end{align*}
 and so,
   $$
      d - I(v) + o_n(1) \ge \frac{1}{p_+} \int_{\mathbb{R}^{N}} (| \nabla w_n |^{ p(x) } + V(x) | w_n |^{ p(x) } ) - \frac{1}{p^{ \ast }_-} \int_{\mathbb{R}^{N}} | w_n |^{ p^{ \ast }(x) } .
   $$
Taking the limit of $n \to +\infty$ in the last inequality, we see that
   \begin{equation} \label{relation3}
      d - I(v) \ge \frac{1}{p_+} L - \frac{1}{p^{ \ast }_-} L = \gamma L.
   \end{equation}
In this moment, it is very important to recall that
$$
    I(v) \geq -M \left( \lambda^{ \Theta_- } + \lambda^{ \Theta_+ } \right).
 $$
 Then, by the hypothesis on $ d $,
\begin{equation} \label{EST1}
d - I(v) < c_\mu.
\end{equation}
On the other hand, since $ \mu \geq \mu_0 $, the last inequality combined with (\ref{ESTIMATIVAPASSO}) leads to
   $$
      d - I(v) < \frac{1}{2K^{ p_+ }} \nu.
   $$
Using this information, we get
$$
\|w_n\| \leq \frac{1}{K} < 1 \,\,\, \forall \,\, n \geq n_0
$$
and so,
$$
| w_n |_{p^{*}(x)} \leq 1 \,\,\, \forall \,\, n \geq n_0.
$$
The above inequalities show that
$$
      \left( \int_{\mathbb{R}^{N}} | w_n |^{ p^{ \ast }(x) }  \right)^{ 1/p^{ \ast }_- } \le K \left(  \int_{\mathbb{R}^{N}} | \nabla w_n |^{ p(x) } + V(x) | v_n |^{ p(x) } \right)^{ 1/p_+ }, \,\,\, \forall n \ge n_0.
$$
   Taking the limit of $n \to +\infty$, we derive
   $$
      L^{ 1/p^{ \ast }_- } \le K L^{ 1/p_+}
   $$
 Supposing by contradiction that $ L > 0 $, we obtain
   \begin{equation} \label{relation4}
      L \ge \left( \frac{1}{K} \right)^{ \frac{1}{\gamma} }.
   \end{equation}
Combining \eqref{relation3} with \eqref{relation4}, it follows that
   $$
      d - I(v) \ge \gamma L \ge \gamma \left( \frac{1}{K} \right)^{ \frac{1}{\gamma} },
   $$
   which is a contradiction, once that \eqref{EST1} and \eqref{ESTIMATIVAPASSO} imply that
   $$
      d - I(v) < \gamma \left( \frac{1}{K} \right)^{ \frac{1}{\gamma} }.
   $$
   Thereby, $ L = 0 $. \qed

\section{Existence of solution with positive energy}

In this section, we will show the existence of a solution via Mountain Pass Theorem. Our first lemma establishes that $I$ verifies the mountain pass geometry.

\begin{lemma} \label{MPG}
   For each $ \mu > 0 $, there exists $ \lambda_1 = \lambda_1(\mu) > 0 $ such that $ I $ satisfies the mountain pass geometry, if $ \lambda \in ( 0, \lambda_1 ) $.
\end{lemma}

\noindent {\bf Proof.}   First of all, we observe that
   \begin{multline*}
      I( u ) \ge \frac{1}{p_+} \int_{\mathbb{R}^{N}} \left( \left| \nabla u \right|^{ p(x) } + V(x) | u |^{ p(x) } \right) - \frac{\lambda}{r_-} \int_{\mathbb{R}^{N}} h(x) | u |^{ r(x) }  \\
      - \frac{\mu}{q_-} \int_{\mathbb{R}^{N}} | u |^{ q(x) } - \frac{1}{p^{ \ast }_-} \int_{\mathbb{R}^{N}} | u |^{ p^{ \ast }(x) } \,\,\,\,\, \forall u \in W^{ 1,p(x) } \big( \mathbb{R}^{N} \big).
   \end{multline*}
By Sobolev's embedding, there are $ C_1, C_2 > 0 $ such that
$$
| u |_{ q(x) } \le C_1 \| u \|  \,\,\, \mbox{and} \,\,\,  | u |_{ p^{ \ast }(x) } \le C_2 \| u \|, \,\,\, \forall u \in W^{ 1,p(x) } \big( \mathbb{R}^{N} \big).
$$
If we suppose that
   $$
      \| u \| < m = \min \left\{ 1, \frac{1}{C_1}, \frac{1}{C_2} \right\},
   $$
   then
   $$
      \| u \| < 1, | u |_{ q(x) } < 1 \,\,\,  \text{and} \,\,\, | u |_{ p^{ \ast }(x) } < 1.
   $$
The above inequalities yield
  $$
      I( u ) \ge \frac{1}{p_+} \| u \|^{ p_+ } - \lambda C_3 | h |_{ \Theta(x) } \| u \|^{ r_- } - \mu C_4 \| u \|^{ q_- } - C_5 \| u \|^{ p^{ \ast }_- },  \,\, \text{if} \ \| u \| < m.
$$
   Since $ p_+ < q_-, p^{ \ast }_- $, we can choose $ R=R(\mu) \in (0,m)$  such that
   $$
      \frac{1}{p_+} R^{ p_+ } - \mu C_4 R^{ q_- } - C_5 R^{ p^{ \ast }_- } \ge \frac{1}{2p_+} R^{ p_+ }.
   $$
   So, if $ \| u \| = R $,
   $$
      I( u ) \ge \frac{1}{2p_+} R^{ p_+ } - \lambda C_3 | h |_{ \Theta(x) } R^{ r_- }.
   $$
Now, we choose $ \lambda_1=\lambda_1(\mu) > 0 $ such that
   $$
      \frac{1}{2p_+} R^{ p_+ } - \lambda_1 C_3 | h |_{ \Theta(x) } R^{ r_- } = \beta > 0.
   $$
Consequently, if $ \lambda \in ( 0 , \lambda_1 ) $, we have that
   $$
      I( u ) \ge \beta, \ \text{for} \ \| u \| = R,
   $$
showing that the  first geometry is satisfied. For the second geometry, we fix $ u \in W^{ 1, p(x) } \big( \mathbb{R}^{N} \big) $ with $ u^+ \ne 0 $. Then, for $ t > 1 $,
   \begin{multline*}
      I( tu ) \le t^{ p_+ } \int_{\mathbb{R}^{N}} \frac{1}{p(x)} \left( \left| \nabla u \right|^{ p(x) } + V(x) | u |^{ p(x) } \right)  - t^{ r_- } \int_{\mathbb{R}^{N}} \frac{h(x)}{r(x)} \left( u^+ \right)^{ r(x) }  \\
      - t^{ q_- } \int_{\mathbb{R}^{N}} \frac{1}{q(x)} \left( u^+ \right)^{ q(x) } - t^{ p^{ \ast }_- } \int_{\mathbb{R}^{N}} \frac{1}{p^{ \ast }(x)} \left( u^+ \right)^{ p^{ \ast }(x) } ,
   \end{multline*}
from where it follows that
   $$
      \lim_{ t \to \infty } I( tu ) = - \infty.
   $$
From this, we observe that the second geometry follows choosing  $ e = t_0u $ with $ t_0 > \frac{R}{\| u \|} $ and $ I( t_0u ) \le 0 $.

\begin{lemma} \label{MPL}
   For each $ \mu \ge \mu_0 $, there exists $ 0 < \lambda_2=\lambda_2(\mu) \le \lambda_1 $, with $ \lambda_1 $ given in Lemma \ref{MPG}, such that the mountain pass level $ c $ of  $ I $ satisfies
   $$
       c < c_\mu - M \left( \lambda^{ \Theta_- } + \lambda^{ \Theta_+ } \right),
   $$
for all  $ \lambda \in (0, \lambda_2) $.
\end{lemma}

\noindent {\bf Proof.}
  For each $ \mu \ge \mu_0 $, we know that there is $ \Psi \in W^{ 1,p(x) } \big( \mathbb{R}^{N} \big) \setminus \left\{ 0 \right\} $ with $ \Psi \ge 0 $ such that
   $$
      I_\mu( \Psi ) = c_\mu \ \text{and} \ I'_\mu( \Psi ) = 0.
   $$
In what follows, fix $ \delta_1 > 0 $ such that
   $$
      c_\mu - M \left( \lambda^{ \Theta_- } + \lambda^{ \Theta_+ } \right)  > \frac{c_\mu}{2}, \, \forall \lambda \in (0, \delta_1).
   $$
Since for $ t > 0 $ sufficiently small
$$
I( t\Psi ) \leq t^{ p_- } \int_{\mathbb{R}^{N}} \frac{1}{p(x)} \left( \left| \nabla u \right|^{ p(x) } + V(x) | u |^{ p(x) } \right),
$$
there is $ t_0 > 0 $, which is independent of $\mu$ and $\lambda$, such that
   $$
      I( t\Psi ) \le \frac{c_\mu}{2}, \,\,\, \forall t \in [0, t_0].
   $$
Therefore, for each $\lambda \in (0, \delta_1)$,
   $$
      I( t\Psi ) \le \frac{c_\mu}{2} <  c_\mu - M \left( \lambda^{ \Theta_- } + \lambda^{ \Theta_+ } \right), \,\,\,\,  \forall t \in [0 , t_0].
   $$
   On the other hand, using the fact that $ \Psi \ge 0 $,  we have that
   $$
      I( t\Psi ) = I_\mu( t\Psi ) - \lambda \int_{\mathbb{R}^{N}} \frac{h(x)}{r(x) } ( t\Psi )^{ r(x) } \,\,\, \mbox{for} \,\,  t \ge 0,
   $$
from where it follows that
 $$
      I( t\Psi ) \le c_\mu - \lambda \min \left\{ t^{ r_- }, t^{ r_+ } \right\} \int_{\mathbb{R}^{N}} \frac{h(x)}{r(x) } \Psi^{ r(x) } .
$$
   In particular, for $ t \ge t_0 $,
   $$
      I( t\Psi) \le c_\mu - \lambda \min \left\{ t_0^{ r_- }, t_0^{ r_+ } \right\} \int_{\mathbb{R}^{N}} \frac{h(x)}{r(x) } \Psi^{ r(x) }.
   $$
Fixing $ \delta_2 > 0 $ such that
   $$
      \lambda^{ \Theta_- - 1} + \lambda^{ \Theta_+ - 1} < \frac{\min \left\{ t_0^{ r_- }, t_0^{ r_+ } \right\}}{M}  \int_{\mathbb{R}^{N}} \frac{h(x)}{r(x) } \Psi^{ r(x) } , \,\, \forall \lambda \in (0,\delta_2),
   $$
we have that
$$
     \sup_{t \geq t_0} I( t\Psi ) < c_\mu - M \left( \lambda^{ \Theta_- } + \lambda^{ \Theta_+ } \right), \,\,\,  \text{if} \ \lambda \in (0, \delta_2).
$$
Setting $ \lambda_2 = \min \left\{ \lambda_1, \delta_1, \delta_2 \right\} $, we obtain by the previous estimates,
   $$
      \sup_{ t \ge 0 } I( t\Psi ) < c_\mu - M \left( \lambda^{ \Theta_- } + \lambda^{ \Theta_+ } \right)  \, \forall \lambda \in (0, \lambda_2).
   $$
Once that
$$
c \leq \sup_{ t \ge 0 } I( t\Psi ),
$$
for $\lambda \in (0, \lambda_2)$, it follows that
   $$
      c < c_\mu - M \left( \lambda^{ \Theta_- } + \lambda^{ \Theta_+ } \right),
   $$
finishing the proof of the lemma. \qed

\begin{theorem} \label{T2}
   For each  $ \mu \ge \mu_0 $, there exists $ \lambda^{ \star }=\lambda^{\star}(\mu) > 0 $ such that problem $ ( P ) $ has a solution with positive energy, for all $ \lambda \in ( 0, \lambda^{ \star } ) $.
\end{theorem}

\noindent {\bf Proof.}   Since $ \mu \ge \mu_0 $, by Lemma \ref{compactness}, the functional  $ I $ verifies the $ (PS)_d $  condition for
   $$
      d < c_\mu - M \left( \lambda^{ \Theta_- } + \lambda^{ \Theta_+ } \right).
   $$
In what follows, we fix $ \lambda^{ \star } = \lambda_2,$ where $ \lambda_2 $ was obtained in Lemma \ref{MPL}. From this, if $ \lambda \in (0, \lambda^{ \star } ) $,  by Lemma \ref{MPG}, $ I $ has the mountain pass geometry, and by Lemma \ref{MPL}, the mountain pass level $ c $
   satisfies
   $$
      0<c < c_\mu - M \left( \lambda^{ \Theta_- } + \lambda^{ \Theta_+ } \right).
   $$
Thereby, $I$ satisfies the $(PS)_c$ condition, and so,  there exists $ \Psi_1 \in W^{ 1,p(x) }(\mathbb{R}^N)$ such that
$$
I'(\Psi_1)=0 \,\,\, \mbox{and} \,\,\, I(\Psi_1)=c >0
$$
showing that $ \Psi_1 $ is a  nontrivial solution for $ (P) $ with positive energy . \qed

\vspace{0.5 cm}

\section{Existence of solution with negative energy}

In this section we will show the existence of a solution with negative energy by using Ekeland's Variational Principle.

\begin{lemma} \label{negative infimum}
   $ I $ is bounded below in $ \overline{ B}_R(0)  $, where $ R > 0 $ is given by Lemma \ref{MPG}. Moreover,
   $$
     J = \inf_{ u \in \overline{ B}_R(0) }  I(u) < 0.
   $$
\end{lemma}
\noindent {\bf Proof.} \, If $ u \in \overline{ B}_R(0) $, then $ \| u \| < 1 $. Arguing  like in the proof of  Lemma \ref{MPG}, we obtain
   \begin{align*}
      | I( u ) | & \le \frac{1}{p_-} \int_{\mathbb{R}^{N}} \left( \left| \nabla u \right|^{ p(x) } + V(x) | u |^{ p(x) } \right)  + \frac{\lambda}{r_-} \int_{\mathbb{R}^{N}} h(x) | u |^{ r(x) }\\
      & \phantom{ \int_{\mathbb{R}^{N}} } + \frac{\mu}{q_-} \int_{\mathbb{R}^{N}} | u |^{ q(x) }  + \frac{1}{p^{ \ast }_-} \int_{\mathbb{R}^{N}} | u |^{ p^{ \ast }(x) }  \\
      & \le \frac{1}{p_-} \| u \|^{ p_- } + \lambda C_3 | h |_{ \Theta(x) } \| u \|^{ r_- } + \mu C_4 \| u \|^{ q_- } + C_5 \| u \|^{ p^{ \ast }_- } \\
      & \le \frac{1}{p_-} R^{ p_- } + \lambda C_3 | h |_{ \Theta(x) } R^{ r_- } + \mu C_4 R^{ q_- } + C_5 R^{ p^{ \ast }_- }.
   \end{align*}
From this, $ I $ is bounded from below in $ \overline{ B}_R(0)  $. \\

   Let $ u \in W^{ 1,p(x) } \big( \mathbb{R}^{N} \big) \setminus \left\{ 0 \right\} $ with $ u^+ \ne 0 $ and $ 0 < t < 1 $. Then,
   \begin{multline*}
      I( tu ) \le t^{ p_- } \rho( u ) - \lambda t^{ r_+ } \int_{\mathbb{R}^{N}} \frac{h(x)}{r(x)} \left( u^+ \right)^{ r(x) }  \\
      - \mu t^{ q_+ } \int_{\mathbb{R}^{N}} \frac{1}{q(x)} \left( u^+ \right)^{ q(x) }  - t^{ p^{ \ast }_+ } \int_{\mathbb{R}^{N}} \frac{1}{p^{ \ast }(x)} \left( u^+ \right)^{ p^{ \ast }(x) }.
   \end{multline*}
   Since $ r_+ < p_-, q_+, p^{ \ast }_+ $,
   $$
      I( tu ) < 0, \ \text{for} \ t \approx 0^+,
   $$
  leading to
   $$
     J = \inf_{ u \in \overline{ B}_R(0)  } I(u) < 0.
   $$

\qed

The next result establishes the existence of a $(PS)_J$ sequence for $I$. The main tool used is Ekeland's Variational Principle and the arguments are very similar to those found in \cite{Alves}, this way, its proof will be omitted.

\begin{lemma} \label{Ekeland} For each $\lambda \in (0, \lambda_1)$, where $ \lambda_1 $ is given by Lemma \ref{MPG}, there is a $(PS)_J$ sequence for $I$, that is, there is $(u_n) \subset W^{1,p(x)}(\mathbb{R}^{N})$ satisfying
$$
I(u_n) \to J \,\,\, \mbox{and} \,\,\, I'(u_n) \to 0
$$

\end{lemma}

Now, we are able to prove the existence of a solution with negative energy.

\vspace{0.5 cm}

\begin{theorem} \label{T3}
   For each  $ \mu \ge \mu_0 $, there exists $ \lambda^{ \star \star } > 0 $ such that problem $ ( P ) $ has a solution with negative energy for all $ \lambda \in ( 0, \lambda^{ \star \star } ) $.
\end{theorem}

\noindent {\bf Proof.}   In fact, once that  $ \mu \ge \mu_0 $, by Lemma \ref{compactness} functional $ I $ verifies the $ (PS)_d $ condition for
   $$
      d < c_\mu - M \left( \lambda^{ \Theta_- } + \lambda^{ \Theta_+ } \right).
   $$
In what follows, we choose $ \lambda_3 > 0 $ such that
   $$
      0 < c_\mu - M \left( \lambda^{ \Theta_- } + \lambda^{ \Theta_+ } \right), \, \forall \lambda \in (0, \lambda_3)
   $$
and $ \lambda^{ \star \star } = \min \left\{ \lambda_1, \lambda_3 \right\}. $ \, For each $ \lambda \in ( 0, \lambda^{ \star \star } ) $, it follows from Lemma \ref{Ekeland} that there exists a $ (PS)_J $ sequence $ ( u_n ) $ for $ I $, where
   $$
     J = \inf_{ u \in \overline{ B}_R(0) }  I(u).
   $$
   By Lemma \ref{negative infimum}, we have  $ J < 0 $, then $ I $ verifies the $ (PS)_J $ condition. From this, there exists $ \Psi_2 \in W^{ 1,p(x) } \big( \mathbb{R}^{N} \big) $ such that
   $$
     I'(\Psi_2)=0 \,\,\, \mbox{and} \,\,\, I(\Psi_2)=J<0
   $$
   Hence, $ \Psi_2 $ is a nontrivial solution for $ ( P ) $ with negative energy.    \qed

\section{Final comments}

Regarding to the problem 
  $$
        \begin{cases}
            -\Delta_{p(x)} u + V(x) |u|^{ p(x)-2 }u = \lambda h(x) |u|^{ r(x)-2 }u + \mu |u|^{ q(x)-2 }u + |u|^{ p^{ \ast }(x)-2 }u, \, \mathbb{R}^{N}, \\
              u \in W^{ 1,p(x) } \big( \mathbb{R}^{N} \big),
         \end{cases} \eqno{(P)_* },
   $$
repeating the same arguments used by Azorero \& Alonso \cite{AA}, we can prove that there exists $ \mu^* > 0 $ with the following property: for each $ \mu \ge \mu^* $, there is
	 $ \lambda_\mu > 0 $ such that $ (P)_* $ has infinitely many solutions with negative energy, if $ \lambda \in (0, \lambda_\mu ) $.  This result is obtained using the concept     and properties of genus and working with a truncation of the energy functional corresponding to $ (P)_* $.

\end{document}